\documentclass{ifacconf}

\usepackage{graphicx}      
\usepackage{natbib}        
\usepackage{amsmath}
\usepackage{amssymb}
\usepackage{pgfplots}
\pgfplotsset{compat=1.5}
\begin{document}
\begin{frontmatter}

\title{Nonlinear Model Order Reduction using Diffeomorphic Transformations of a Space-Time Domain\thanksref{footnoteinfo}}

\thanks[footnoteinfo]{The authors acknowledge support by the Deutsche Forschungsgemeinschaft under Germany’s Excellence Strategy EXC 2044 390685587, Mathematics Münster: Dynamics – Geometry – Structure.}

\author[First]{Hendrik Kleikamp} 
\author[First]{Mario Ohlberger} 
\author[First]{Stephan Rave}

\address[First]{Institute for Analysis and Numerics, University of Münster, Germany (e-mail: hendrik.kleikamp@uni-muenster.de).}


\end{frontmatter}

\section{Introduction}
Hyperbolic conservation laws play an important role in many applications. For instance, modeling the behavior of fluids or gases leads to conservation equations for quantities like mass, momentum, or energy. The resulting equations are typically nonlinear and exhibit phenomena like shock formation and transport. Parametrized hyperbolic equations suffer from a highly nonlinear solution manifold that cannot be approximated appropriately by a linear subspace, that is, the solution manifold has a slowly decaying Kolmogorov $N$-width, see for instance \cite{ohlberger2016}. Therefore, methods that rely solely on linear combinations of ansatz-functions are not sufficient to achieve suitable reduced models. Furthermore, the formation and interaction of shocks is an additional difficulty when dealing with model order reduction for conservation laws.

In this contribution, we describe a new nonlinear model order reduction technique. The nonlinearity of the approach stems from the exponential map applied to vector fields in Euclidean space. We will thus identify diffeomorphisms, resulting from the application of the exponential map, with vector fields. Afterwards, we use standard ideas from linear model order reduction to compute a subspace of the space of vector fields. During the online phase, for a given parameter, elements from this subspace are computed and the exponential map is used to determine the corresponding diffeomorphism. This transformation is subsequently applied to a (fixed) space-time solution snapshot to obtain the approximate solution for the new parameter.

\section{Basics from differential geometry and image registration}

In this section, we give a brief overview of important notions from differential geometry and introduce the basic concepts of image registration via geodesic shooting.

\subsection{Differential geometry and Lie groups}

In the intersection of differential geometry and group theory, one considers so called \emph{Lie groups}, which are groups such that group multiplication and inversion are smooth. Directly connected to the concept of Lie groups is the notion of \emph{Lie algebras}. For a Lie group $G$, the corresponding Lie algebra $\mathfrak{g}$ is defined as the tangent space to the manifold $G$ at the identity element. The \emph{exponential map} $\exp\colon\mathfrak{g}\to G$ describes, for an element $v\in\mathfrak{g}$ of the Lie algebra, the end point of a shortest path (a \emph{geodesic}) that starts at the identity in $G$ in the direction given by $v$.

\subsection{Image registration and geodesic shooting}

The field of image registration has its origins in the analysis of medical image data. Given two images $u_0,u_1\colon\Omega\to\mathbb{R}^d$, treated as functions on a fixed domain $\Omega\subseteq\mathbb{R}^d$, the aim of image registration is to find a transformation $\phi\colon\Omega\to\Omega$ such that $u_0\circ\phi^{-1}\approx u_1$. There exist several choices for the class of transformation to employ. A quite general approach uses the group of diffeomorphisms of the domain $\Omega$. Since the group of diffeomorphisms also forms a Lie group, with the vector space of smooth vector fields on $\Omega$ being the corresponding Lie algebra, we parametrize a diffeomorphism by a single vector field to which we apply the exponential map to regain the corresponding diffeomorphism. This idea is used in the \emph{geodesic shooting} algorithm developed in \cite{miller2006}. In the aforementioned work, the Euler-Poincaré equations for the evolution along a geodesic in the diffeomorphism group are described. Together with the corresponding adjoint equations, it is possible to formulate a gradient descent algorithm for an energy functional of the form
\[
	E_{u_0\to u_1}(v_0) := \lVert v_0\rVert_V^2 + \frac{1}{\sigma^2}\lVert u_0\circ\phi_1^{-1}-u_1\rVert_{L^2(\Omega)}^2,
\]
where $v_0\colon\Omega\to\mathbb{R}^d$ denotes a vector field, $\lVert\cdot\rVert_V$ is a suitable norm on the space of vector fields, $\sigma>0$ is a weighting parameter, and the diffeomorphism $\phi_1\colon\Omega\to\Omega$ is given as the solution at the final time $t=1$ of the equation
\[
	\frac{d\phi_t}{dt} = v_t\circ\phi_t,
\]
where $v_t$, for $t\in[0,1]$, solves the Euler-Poincaré equation for the initial vector field $v_0$. It then holds $\phi_1=\exp(v_0)$.

\section{Nonlinear approximation scheme}

Before describing the algorithm in detail, we introduce some more notation: We denote by $\mathcal{P}\subset\mathbb{R}^p$ for some $p\in\mathbb{N}$ the parameter space. Moreover, the space-time solution of the equation under consideration for the parameter $\mu\in\mathcal{P}$ is denoted by $u(\mu)\colon\Omega\to\mathbb{R}$. The domain $\Omega$ is a subset of $\mathbb{R}^n\times\mathbb{R}^+$, where $n\in\mathbb{N}$ is the space dimension. It therefore holds $d=n+1$ for the dimension of the domain we deform.

For simplicity, we restrict our attention to a single "topology" of the solution, which is assumed to be independent of the parameter $\mu\in\mathcal{P}$. We consider, for instance, only solutions with a single shock or solutions with two merging shocks. This means that we assume that it is possible to transform solutions into each other by means of diffeomorphic transformations of the underlying domain $\Omega$.

The idea of using Lie groups together with their corresponding Lie algebra for model order reduction of hyperbolic equations was introduced in \cite{ohlberger2013}, where finite-dimensional groups acting only on the spatial domain were considered, for instance the translation group. Here, we use the infinite-dimensional diffeomorphism group on the space-time domain $\Omega$, such that shock formation and interaction are already included in the ansatz-functions.

\subsection{Offline algorithm}

During the offline phase, we first of all choose a reference parameter $\mu_{\text{ref}}\in\mathcal{P}$ and compute the related full-order space-time reference solution $u(\mu_\text{ref})$. Afterwards, we select training parameters $\mu_1,\dots,\mu_M\in\mathcal{P}$ and compute the solution snapshots $u(\mu_1),\dots,u(\mu_M)$. We do not detail the exact solution algorithm for the full-order computations, the only requirement we impose is that the solution data can be treated as a function on $\Omega$, such that we can apply the geodesic shooting algorithm for image registration. Subsequently, vector fields $v(\mu_1),\dots,v(\mu_M)\colon\Omega\to\mathbb{R}^d$ are computed, using the geodesic shooting algorithm, such that they minimize $E_{u(\mu_{\text{ref}})\to u(\mu_1)},\dots,E_{u(\mu_{\text{ref}})\to u(\mu_M)}$. The set of vector fields $v(\mu_1),\dots,v(\mu_M)$ is now reduced using \emph{proper orthogonal decomposition}, similar to the procedure described in \cite{wang2019a}. This step results in an orthogonal matrix $V_N$, whose columns span an $N$-dimensional subspace of the space of vector fields. Finally, an artificial neural network $\Phi\colon\mathcal{P}\to\mathbb{R}^N$ is trained to approximate the mapping $\pi\colon\mathcal{P}\to\mathbb{R}^N$ defined as $\pi(\mu)=V_N^\top v(\mu)$. The function $\pi$ maps a parameter $\mu\in\mathcal{P}$ to the coefficients (with respect to the basis $V_N$) of the orthogonal projection of the optimal vector field $v(\mu)$ onto the subspace $\operatorname{ran}(V_N)$.

\subsection{Online algorithm}

Given a new parameter $\mu\in\mathcal{P}$, a forward pass through the neural network is performed to obtain the approximate coefficients $\Phi(\mu)\approx\pi(\mu)$. Next, the vector field $v_N(\mu)=V_N\Phi(\mu)$ is computed. By applying the exponential map $\exp$ to $v_N(\mu)$, we derive the diffeomorphism $\phi_N(\mu)=\exp(v_N(\mu))$. The approximate solution for the parameter $\mu$ is now given as $u_N(\mu)=u(\mu_{\text{ref}})\circ\phi_N^{-1}(\mu)$.

\section{Example}

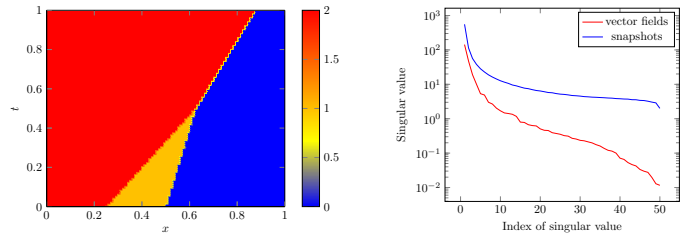
\begin{figure}[h]
	\begin{center}
		\includegraphics[width=.5\linewidth]{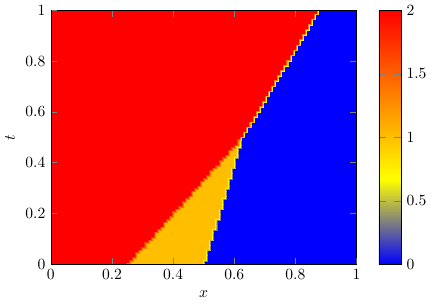}
		\hfill
		\begin{tikzpicture}[scale=0.45]
			\begin{semilogyaxis}[xlabel=Index of singular value, ylabel=Singular value]
				\addplot[color=red] table[x expr=\coordindex+1, y index=0] {singular_values_vector_fields.txt};
				\addlegendentry{vector fields}
				\addplot[color=blue] table[x expr=\coordindex+1, y index=0] {singular_values_snapshots.txt};
				\addlegendentry{snapshots}
			\end{semilogyaxis}
		\end{tikzpicture}
		\caption{Sample solution to Burgers' equation and singular values of vector fields (red) and snapshots (blue)}
		\label{fig:space-time-solution}
	\end{center}
\end{figure}

We present the decay of the singular values of the computed vector fields for a Burgers' equation with two merging shocks. The equation of interest reads
\[
	\partial_t u + \frac{\mu}{2}u\,\partial_x u = 0, \quad u(x,0)=\begin{cases}2, & \text{if }x\leq 0.25,\\ 1, & \text{if }0.25<x<0.5, \\ 0, & \text{otherwise},\end{cases}
\]
where $\mu\in[0.25,1]=:\mathcal{P}$, and $(x,t)\in[0,1]^2=:\Omega$. An example of a space-time solution for this equation for $\mu=1/2$ is given in the left part of Fig.~\ref{fig:space-time-solution}. Starting with the reference parameter $\mu_{\text{ref}}=0.25$, we performed registration onto $50$ snapshots for parameters uniformly selected from $\mathcal{P}$. The singular values of the vector fields together with the singular values of the space-time snapshots themselves are presented in the right part of Fig.~\ref{fig:space-time-solution}. The maximum relative $L_2$-error of the transformed snapshots with respect to the exact solutions is roughly $6\%$. The plots show that the singular values of the vector fields decay much faster (even exponentially) than those of the snapshots, which means that the vector fields can be approximated more efficiently by a linear subspace than the snapshots.

\section{Conclusion}

In this work we describe a new approach for nonlinear model order reduction for parametrized hyperbolic equations. Future research in this direction will be concerned with the computation of the reduced coefficients for the vector fields by solving a residual-minimization problem. Furthermore, the treatment of different solution topologies would make the algorithm more flexible.

\bibliography{ifacconf}

\end{document}